\documentclass[a4paper]{amsart}
\usepackage{amssymb}

\usepackage{graphicx}
\usepackage[matrix,arrow,curve,tips]{xy}
           \SelectTips{cm}{}

\newtheorem{dfn}{Definition}[section]
\newtheorem{prop}[dfn]{Proposition}
\newtheorem{theo}[dfn]{Theorem}
\newtheorem{cor}[dfn]{Corollary}
\newtheorem{lem}[dfn]{Lemma}
\newtheorem{rem}[dfn]{Remark}
\newtheorem{ex}[dfn]{Example}
\newtheorem{exs}[dfn]{Examples}

\newcommand{\ra}{\rightarrow}
\newcommand{\lra}{\longrightarrow}
\newcommand{\oo}{\,\mbox{-}\,}
\newcommand{\com}{\mathbin{{\scriptstyle \circ }}}

\newcommand{\RR}{\mathbb{R}}

\newcommand{\cF}{\mathord{\mathcal{F}}}
\newcommand{\cG}{\mathord{\mathcal{G}}}

\newcommand{\fg}{\mathord{\mathfrak{g}}}
\newcommand{\fh}{\mathord{\mathfrak{h}}}
\newcommand{\fX}{\mathord{\mathfrak{X}}}
\newcommand{\fL}{\mathord{\mathfrak{L}}}

\newcommand{\eC}{\mathord{\mathit{C}^{\infty}}}

\newcommand{\src}{\mathord{\mathrm{s}}}
\newcommand{\trg}{\mathord{\mathrm{t}}}
\newcommand{\uni}{\mathord{\mathrm{u}}}
\newcommand{\anchor}{\mathord{\mathrm{an}}}

\newcommand{\Mon}{\mathord{\mathrm{Mon}}}
\newcommand{\Hol}{\mathord{\mathrm{Hol}}}

\begin{document}

\title{On the integrability of subalgebroids}

\author{I. Moerdijk}
\address{Mathematical Institute, Utrecht University,
         P.O. Box 80.010, 3508 TA Utrecht, The Netherlands}
\email{moerdijk@math.uu.nl}

\author{J. Mr\v{c}un}
\address{Department of Mathematics, University of Ljubljana,
         Jadranska 19, 1000 Ljubljana, Slovenia}
\email{janez.mrcun@fmf.uni-lj.si}
\thanks{This work was supported in part by the
        Dutch Science Foundation (NWO) and
        the Slovenian Ministry of Science (M\v{S}Z\v{S} grant J1-3148)}


\begin{abstract}
Let $G$ be a Lie groupoid with Lie algebroid $\fg$.
It is known that, unlike in the case of Lie groups, not
every subalgebroid of $\fg$ can be integrated by a subgroupoid of $G$.
In this paper we study conditions on the invariant foliation
defined by a given subalgebroid under which such an integration
is possible. We also consider the problem of integrability
by closed subgroupoids, and we give conditions under which
the closure of a subgroupoid is again a subgroupoid.
\end{abstract}

\maketitle

\section*{Introduction} \label{sec:intro}

The basic theory of Lie groups and Lie algebras, contained in any
course on the subject, establishes an almost perfect correspondence
between finite-dimensional Lie algebras and Lie groups.
In particular, any
abstract Lie algebra can be `integrated' to a Lie group, and this group
is uniquely determined when one requires it to be simply connected.
Morphisms from a simply connected Lie group to another one are in
bijective correspondence with morphisms between their Lie algebras. If
$\fg$ is the Lie algebra of a Lie group $G$, subalgebras of $\fg$
correspond to
immersed subgroups of $G$. The closure of such an immersed Lie subgroup
is again one, and an immersed subgroup is embedded if and
only if it is closed.

The basic Lie theory becomes much more involved (and interesting) in
the wider context where one allows local symmetries, i.e.\ in the
context of so-called Lie groupoids and their algebroids. These objects
arise naturally in the theory of actions of Lie groups or algebras on
manifolds, in foliation theory, in Poisson geometry, in areas related
to mathematical physics such as the theory of (gerbes on) orbifolds
and quantization theory, and in many other situations
\cite{BursztynWeinstein2004,CannasdasilvaWeinstein1999,
      CattaneoFelder2000,Connes1994,
      CrainicFernandes2003,Landsman1998,Mackenzie1987,Pradines1967}.
Like for groups, any Lie groupoid $G$ has an infinitesimal counterpart, its Lie
algebroid $\fg=\fL(G)$. Again, morphisms of Lie groupoids correspond
bijectively to morphisms between their Lie algebroids under suitable
connectivity assumptions \cite{MackenzieXu2000,MoerdijkMrcun2002},
but this is about as far as the
simple picture of groups extends to groupoids. The main difference lies
in the fact that not every Lie algebroid is integrable, i.e.\  arises as
the Lie algebroid of a Lie groupoid. This was first noticed by Almeida
and Molino \cite{AlmeidaMolino1985,Molino1988},
who proved that the developability of a certain
kind of foliation was obstructed by the integrability of a
naturally associated algebroid, a result which immediately led to
natural counterexamples to integrability. Since then, our understanding
of integrability of algebroids has advanced considerably, and we refer
the reader to \cite{CrainicFernandes2003}
and the references therein for a recent state of the art.

This note and its sequel \cite{MoerdijkMrcun2004b}
are concerned with the integrability
problem for  subalgebroids. To explain our results,
consider a Lie groupoid $G$ with algebroid $\fg$.
It is known that any
subalgebroid $\fh\subset\fg$
of such an integrable algebroid $\fg$ is again integrable,
by a Lie groupoid $H$ which admits an immersion $H\ra G$.
However, unlike
the simple case of groups and algebras, there may not exist
a {\em subgroupoid} of $G$ integrating $\fh$, or more precisely,
there may not exist an {\em injective} such immersion $H\ra G$
integrating $\fh\subset\fg$. And
even if there is such an injective immersion,
the closure of its
image may not be a Lie subgroupoid of $G$. In this paper, we will use
foliation theory to investigate under which conditions results like in
the case of groups hold in the context of groupoids.
The starting point is that any
subalgebroid $\fh\subset\fg$
gives rise to a right invariant foliation on $G$.
We prove that
$\fh\subset\fg$ can be integrated to an injective immersion
$H\ra G$ if and only if this
foliation has trivial holonomy (a condition automatically fulfilled if
$G$ is a group).
We prove also
that if the invariant foliation associated to $\fh$ is
transversely complete (again, automatically true in the group case),
then $\fh$ integrates to an injective submersion
$\iota\!:H\ra G$, with the property that the closure of its image is a Lie
subgroupoid $\bar{H}\subset G$. Moreover, the image of $\iota$
itself is closed if and only if $\iota$ is an embedding.

These results will be used in a sequel \cite{MoerdijkMrcun2004b}
to this paper, where we will consider the question
when a given Lie subalgebroid $\fh$ of $\fg$ can be integrated
to a {\em closed} subgroupoid of a
possibly larger groupoid $\tilde{G}$
integrating $\fg$. Under some assumptions,
we will prove that this is the case
if and only if a specific algebroid, naturally associated to $\fh$
and $G$, is integrable. This result extends
the classical Almeida-Molino theorem referred to above.

\section{Preliminaries} \label{sec:pl}

For the convenience of the reader, and to fix the notation,
we begin by summarizing some basic definitions.
For detailed exposition and many examples, the reader
might wish to consult one of the books
\cite{CannasdasilvaWeinstein1999,Mackenzie1987,MoerdijkMrcun2003,Molino1988}
and references cited there.

\subsection{Lie algebroids} \label{seusec:la}

Let $M$ be a smooth manifold. We recall that a {\em Lie algebroid}
over $M$ is a (real) vector bundle $\fg$ over $M$, equipped with
a vector bundle map $\anchor\!:\fg\ra T(M)$ (called the {\em anchor})
and a Lie algebra structure on the vector space $\Gamma(\fg)$ of
sections of $\fg$. This structure has to satisfy two axioms:
On sections, the map $\anchor\!:\Gamma(\fg)\ra\fX(M)$ should be
a homomorphism of Lie algebras; furthermore, 
the Leibniz rule
$$ [X,fY]=f[X,Y]+\anchor(X)(f)Y $$
should hold for any $X,Y\in\Gamma(\fg)$ and any $f\in\eC(M)$.
A {\em morphism} $\fg\ra\fh$ of Lie algebroids over $M$ is
a map of vector bundles which preserves the structure.
(There is also a more involved notion of morphism between
algebroids over different manifolds,
see \cite{HigginsMackenzie1990}.)
Such Lie algebroids arise in many contexts: Any foliated manifold
$(M,\cF)$ can be viewed as an algebroid, whose anchor
is the inclusion $\cF\ra T(M)$. Any Poisson manifold $P$
can be viewed as an algebroid, whose anchor $T^{\ast}(P)\ra T(P)$
is defined in terms of the Poisson bracket on the functions
by $\anchor(df)(g)=-\{f,g\}$. Any infinitesimal action of a finite
dimensional Lie algebra $\fg$ on a manifold $M$ defines a
Lie algebroid $\fg\ltimes M$, which is
the trivial vector bundle $\fg\times M\ra M$
with the anchor given by the infinitesimal action.

\subsection{Lie groupoids and integrability} \label{sebsec:lgi}

A Lie groupoid $G$ over a manifold $M$
will be denoted
$$
\xymatrix{
G \ar@<2pt>[r]^{\src} \ar@<-2pt>[r]_{\trg} & M \ar[r]^{\uni} & G
}\;,
$$
where $\src$ and $\trg$ are the source and target maps,
and $\uni$ maps each point $x\in M$ to the unit arrow $1_{x}$.
The multiplication or composition of two arrows $g\!:x\ra y$
and $h\!:y\ra z$ is denoted by $hg\!:x\ra z$. The maps $\src$
and $\trg$ are required to be submersions. We will generally
assume that $M$ and the fibers of $\src$
(and hence also of $\trg$)
are Hausdorff manifolds, but examples force us to allow for
the possibility that $G$ is non-Hausdorff.
However, we assume that the fibers of the source map
are Hausdorff.

For such a groupoid $G$, the bundle $T^{\src}(G)$
of source-vertical tangent vectors pulls back along $\uni\!:M\ra G$
to a vector bundle $\uni^{\ast}(T^{\src}(G))$ over $M$,
which has the structure of a Lie algebroid. This algebroid is
called the Lie algebroid of the groupoid $G$, and denoted by
$\fL(G)$, or sometimes simply by $\fg$. Lie algebroids which arise
in this way are called {\em integrable}.

A Lie groupoid $G$ is called {\em source-connected}
if $\src\!:G\ra M$ has connected fibers. Every Lie groupoid $G$
contains an open subgroupoid whose
source-fibers are the connected components of the source-fibers
of $G$ containing the units. This groupoid has the same
Lie algebroid as $G$. For this reason we can restrict
to source-connected Lie groupoids:
All Lie groupoids considered in this paper
will be assumed source-connected.

For any (source-connected) Lie groupoid $G$ there exists a unique
(up to unique isomorphism) cover $p\!:\tilde{G}\ra G$ by a
source-connected Lie groupoid
$\tilde{G}$ whose source-fibers are simply connected
(we say that $\tilde{G}$ is {\em source-simply connected}).
The source-fiber of $\tilde{G}$ is the universal cover
of the corresponding source-fiber of $G$, and $\tilde{G}$ again
determines the same Lie algebroid as $G$, i.e.\
the induced morphism of Lie algebroids
$\fL(p)\!:\fL(\tilde{G})\ra\fL(G)$ is an isomorphism.

As an example, and for later reference, we mention that
any foliation $\cF$ of $M$, viewed as a Lie algebroid, is
always integrable. It is the Lie algebroid of the
{\em holonomy groupoid} $\Hol(M,\cF)$ of the foliation, and
also of its source-simply connected cover, which is known as the
{\em monodromy groupoid} $\Mon(M,\cF)$ of the foliation.
The source-fibers of the latter groupoid are the universal covers
of the leaves of $\cF$.
(For details, see \cite{MoerdijkMrcun2003}.)

\section{Subalgebroids and invariant foliations} \label{sec:saif}

Throughout this section $G$ denotes a fixed (source-connected)
Lie groupoid over $M$,
and $\fg$ denotes its Lie algebroid. Recall that
all Lie groupoids considered are assumed source-connected.

A {\em subalgebroid} of $\fg$ is a vector subbundle
$\fh\subset\fg$ for which $\Gamma(\fh)\subset\Gamma(\fg)$ is
a Lie subalgebra.
This makes $\fh$ into a Lie algebroid over $M$,
and the inclusion $\fh\ra\fg$ is a morphism of Lie algebroids.
(In this paper, we shall only consider
subalgebroids over the same base.)
We recall from \cite{MoerdijkMrcun2002} that any subalgebroid
of an integrable Lie algebroid is itself integrable,
and in particular $\fh$ integrates to a Lie groupoid $H$, while
the morphism $\fh\ra\fg$ integrates to an immersion
$H\ra G$ of Lie groupoids. Our first goal is to describe
the range of immersed groupoids which arise as integrals
of the given subalgebroid $\fh$.

Let $\cF(\src)=T^{\src}(G)$ be the foliation of $G$ by
the fibers of the source map. There is a canonical isomorphism
$$
\trg^{\ast}(\fg)\stackrel{\cong}{\lra} \cF(\src)
$$
of vector bundles over $G$, and $\fh\subset\fg$ defines another foliation
$\cF(\fh)\cong \trg^{\ast}(\fh)$ which refines $\cF(\src)$,
i.e.\ $\cF(\fh)\subset\cF(\src)$.
The groupoid $G$ acts on itself by right multiplication
along the source map, and $\cF(\fh)$ is invariant
under this action. Conversely, any such right invariant
foliation $\cF\subset\cF(\src)$ is of the form
$\cF(\fh)$ for a unique subalgebroid $\fh$,
defined in terms of $\cF$ and the unit section $\uni\!:M\ra G$
by $\fh=\uni^{\ast}(\cF)$. For later reference,
we record this correspondence.

\begin{lem} \label{lem:saif.3}
Let $G$ be a
Lie groupoid with Lie algebroid $\fg$.
Subalgebroids of $\fg$ correspond bijectively to the right
invariant foliations $\cF\subset\cF(\src)$ of $G$.
\end{lem}

By the invariance of such a foliation $\cF(\fh)$, the
groupoid $G$ is expected to act naturally on structures
associated to $\cF(\fh)$. In particular:

\begin{lem} \label{lem:saif.5}
Let $G$ be a
Lie groupoid with Lie algebroid $\fg$,
and let $\fh$ be a subalgebroid of $\fg$. Then
$G$ naturally acts on the monodromy and holonomy
groupoids of the associated right invariant foliation
$\cF(\fh)$ of $G$.
\end{lem}

\begin{proof}
The groupoids $\Hol(G,\cF(\fh))$ and $\Mon(G,\cF(\fh))$
are groupoids over the space $G$,
and this space carries and action of
the groupoid $G$ by right
multiplication. An arrow in $\Mon(G,\cF(\fh))$
is the homotopy class $[\alpha]$ of a path
$\alpha$ inside a leaf of $\cF(\fh)$, which is
contained in a source-fiber $\src^{-1}(x)$ of $G$.
An arrow $g\!:y\ra x$ of $G$ acts on $[\alpha]$ in the obvious way
by $[\alpha]g=[\alpha g]$, where
$(\alpha g)(t)=\alpha(t) g$ for any $t\in [0,1]$.
This is a path in a leaf of $\cF(\fh)$ because
$\cF(\fh)$ is right invariant, and this action is obviously
well-defined on homotopy classes. This describes the
action of $G$ on $\Mon(G,\cF(\fh))$.

We only need to check
that this action also descends to an action on
$\Hol(G,\cF(\fh))$.
To see this, choose a local bisection
$\sigma\!:U\ra G$, defined on a neighbourhood
$U$ of $y$, with $\sigma(y)=g$. Now the action
by elements of $\sigma(U)$ gives us a diffeomorphism
$\src^{-1}(\trg(\sigma(U)))\ra \src^{-1}(U)$ which
preserves the foliation $\cF(\fh)$.
This in particular implies that if $\alpha$ is a loop
with trivial holonomy, then so is $\alpha g$,
so the action of $G$ on $\Mon(G,\cF(\fh))$
descends to an action on $\Hol(G,\cF(\fh))$
for which the quotient map
$\Mon(G,\cF(\fh))\ra \Hol(G,\cF(\fh))$
is $G$-equivariant.
\end{proof}

Since the $G$-action on itself is principal,
so are the $G$-actions on $\Mon(G,\cF(\fh))$ and
$\Hol(G,\cF(\fh))$. Therefore
\cite{MoerdijkMrcun2002,MoerdijkMrcun2003}
we can form the quotient Lie groupoids
$$ H_{\mathrm{max}}=\Mon(G,\cF(\fh))/G $$
and
$$ H_{\mathrm{min}}=\Hol(G,\cF(\fh))/G $$
over $M$.

\begin{theo} \label{theo:saif.6}
Let $G$ be a Lie groupoid 
with source-simply connected cover $\tilde{G}$
and with Lie algebroid $\fg$,
and let $\fh$ be a subalgebroid of $\fg$.

(i) The Lie  groupoids
$H_{\mathrm{max}}$
and $H_{\mathrm{min}}$ defined above
both integrate the subalgebroid $\fh$, and fit into 
a natural commutative square
$$
\xymatrix{
H_{\mathrm{max}} \ar[d] \ar[r] & H_{\mathrm{min}} \ar[d] \\
\tilde{G} \ar[r] & G
}
$$
where both vertical maps are immersions integrating the inclusion
$\fh\ra\fg$.

(ii) 
Let $\iota\!:H\ra G$ be any immersion of Lie groupoids
over $M$ which integrates the inclusion $\fh\ra\fg$.
Then there are natural maps
$$ H_{\mathrm{max}} \lra H \lra H_{\mathrm{min}} $$
of Lie groupoids over $M$ which both integrate the identity
on $\fh$, and $\iota$ factors through
$H \ra H_{\mathrm{min}}$ as the canonical immersion
$H_{\mathrm{min}}\ra G$.
\end{theo}

\begin{rem} \rm \label{rem:saif.7}
Motivated by this theorem, the Lie groupoid
$H_{\mathrm{min}}$ will be referred as the
{\em minimal integral} of $\fh$ over $G$.
Analogously, the Lie groupoid
$H_{\mathrm{max}}$ will be called the
{\em maximal integral} of $\fh$.
Note, however, that $H_{\mathrm{max}}$ is in fact
the source-simply connected integral $\tilde{H}$ of $\fh$.
In particular, it is independent of the choice of
$G$ (up to isomorphism).
Since we are assuming that $H$ is source-connected,
the maps in (ii) are surjective submersions.
Moreover, they restrict to covering projections on source-fibers.
\end{rem}

\begin{proof}[Proof of Theorem \ref{theo:saif.6}]
(i)
That $H_{\mathrm{max}}$ integrates $\fh$
follows easily from the observation
\cite[Remark on p. 573]{MoerdijkMrcun2002}
that the quotient map
$\Mon(G,\cF(\fh))\ra H_{\mathrm{max}}$ induces an isomorphism
of source-fibers.
Exactly the same argument applies to the map
$\Hol(G,\cF(\fh))\ra H_{\mathrm{min}}$.
Next, since $\cF(\fh)$ refines the simple foliation
$\cF(\src)$, there is a natural square
$$
\xymatrix{
\Mon(G,\cF(\fh)) \ar[d] \ar[r] & \Hol(G,\cF(\fh)) \ar[d] \\
\Mon(G,\cF(\src)) \ar[r] & \Hol(G,\cF(\src))
}
$$
of groupoids over $G$ and $G$-equivariant morphisms
between them. Factoring out the principal $G$-action,
we obtain the diagram of groupoids over $M$ as in the statement.

(ii)
Suppose that $\iota\!:H\ra G$ is an immersion which integrates
$\fh\ra\fg$. Consider the groupoid $\hat{H}$ over $G$ whose
arrows $g\ra g'$ are arrows $h$ of $H$ with $\iota(h)g=g'$.
Then $\hat{H}$ integrates $\cF(\fh)$ and carries
an obvious principal $G$-action for which $\hat{H}/G=H$.
By \cite[Proposition 1]{CrainicMoerdijk2001}
there are maps
$$\Mon(G,\cF(\fh)) \lra \hat{H} \lra \Hol(G,\cF(\fh)) $$
of foliation groupoids over $G$ which all integrate
the same foliation $\cF(\fh)$. Factoring out
the $G$-action, we obtain the desired maps
of Lie groupoids over $M$.
\end{proof}

\begin{cor} \label{cor:saif.8}
Let $G$ be a 
Lie groupoid with
Lie algebroid $\fg$, and let $\fh$ be a subalgebroid of $\fg$
with minimal integral $H_{\mathrm{min}}$ over $G$.
Then the following conditions are equivalent:
\begin{itemize}
\item [(i)]   The foliation $\cF(\fh)$ has trivial holonomy.
\item [(ii)]  The canonical immersion
              $H_{\mathrm{min}}\ra G$ is injective.
\item [(iii)] The inclusion $\fh\ra\fg$ can be integrated
              to an injective immersion $H\ra G$.
\end{itemize}
\end{cor}

\begin{proof}
Note first that $\cF(\fh)$ has trivial holonomy if
and only if the map 
$$ \Hol(G,\cF(\fh))\lra\Hol(G,\cF(\src)) $$
is injective. Since the $G$-actions on $\Hol(G,\cF(\fh))$
and on $\Hol(G,\cF(\src))$ are principal, this is equivalent
to injectivity of the map
$H_{\mathrm{min}}\ra G$, which integrates
the inclusion $\fh\ra\fg$ by 
Theorem \ref{theo:saif.6} (i).
This shows equivalence between (i) and (ii).

Since (ii) is clearly stronger than (iii),
we only need to show that (iii) implies (ii).
Indeed, if $H\ra G$ is any injective immersion
integrating the inclusion $\fh\ra\fg$, we can
assume that $H$ is
source-connected, and we can use
Theorem \ref{theo:saif.6} (ii) to obtain a map
$H\ra H_{\mathrm{min}}$ which integrates the identity
on $\fh$. Furthermore, since the
composition
$$ H\lra H_{\mathrm{min}}\lra G $$
is injective by assumption, and since
$H\ra\Hol(G,\cF(\fh))$ is a covering projection
on each source-fiber \cite{CrainicMoerdijk2001},
the map $H_{\mathrm{min}}\ra G$ is injective as well.
\end{proof}

\begin{exs} \rm \label{exs:saif.9}
(1)
Let $G$ be a connected Lie group and $\fg$ its Lie algebra of right
invariant vector fields.
Then any subalgebra $\fh\subset\fg$ determines a right invariant
foliation $\cF(\fh)$, whose tangent space at $g\in G$ is
$\cF(\fh)_{g}=dR_{g}(\fh)\subset T_{g}(G)$. As is well known,
$\fh$ is integrated by a connected
Lie group $H$ which occurs as the leaf
of $\cF(\fh)$ through $1\in G$, so obviously $H$ is injectively immersed
in $G$. By Corollary \ref{cor:saif.8} the foliation $\cF(\fh)$
has trivial holonomy, a fact which is well-known in this case.

(2)
Let $\cF$ be a foliation of a manifold $M$.
We can view $\cF$ as a subalgebroid of $T(M)$,
which is the algebroid of the pair groupoid $M\times M$ over $M$.
The source-simply connected cover of the pair groupoid is
the fundamental groupoid $\Pi(M)$ of $M$. The foliation $\cF$
can be integrated by an injectively immersed groupoid
$H\ra M\times M$ if and only if it has trivial holonomy, and
by an injectively immersed groupoid $H'\ra \Pi(M)$ if and only
if the pull-back of $\cF$ to the universal cover $\tilde{M}\ra M$
has trivial holonomy.

(3)
Consider the torus $T^{2}=S^{1}\times S^{1}$, which is
a Lie group with commutative Lie algebra $\RR^{2}$.
Any $\theta\in\RR P^{1}$ is a
one-dimensional subspace of $\RR^{2}$ and hence a
Lie subalgebra, so by (i) it determines a foliation
$\cF(\theta)$ of $T^{2}$ with trivial holonomy
(it is the Kronecker
foliation of the torus with slope $\theta$).
There is a countable dense subset $A$ of $\RR P^{1}$
such that $\cF(\theta)$ has only compact leaves for
$\theta \in A$ and only non-compact dense leaves for
$\theta \not\in A$.

These foliations together define a
foliation $\cF$ of the trivial bundle of Lie groups
$T^{2}\times\RR P^{1}$ over $\RR P^{1}$. If we view
this bundle as a Lie groupoid over $\RR P^{1}$,
then the leaves of $\cF$ are contained in the source fibers
and $\cF$ is right invariant. Therefore it determines a
subalgebroid $\fh$ of the Lie algebroid $\fg=\RR^{2}\times \RR P^{1}$
associated to the groupoid $T^{2}\times\RR P^{1}$.
This subalgebroid is not integrable by a groupoid injectively
immersed in $T^{2}\times\RR P^{1}$. Indeed, the
local Reeb stability theorem implies that the holonomy
of $\cF$ is not trivial.

On the other hand, note that $\fg$ is also integrable by the
source-simply connected groupoid $\RR^{2}\times\RR P^{1}$,
and that the subalgebroid $\fh$ is integrable by a groupoid
injectively immersed in $\RR^{2}\times\RR P^{1}$.
\end{exs}

\section{Subgroupoids of Lie groupoids} \label{sec:sglg}

Let $G$ be a fixed Lie groupoid over a
connected manifold $M$, and write $\fg$ for the Lie algebroid
of $G$.
Recall that all Lie groupoids considered are assumed
source-connected.

A {\em subgroupoid} of $G$ is another
(source-connected) Lie groupoid $H$ over $M$,  equipped with
an injective immersion $\iota\!:H\ra G$, which is also
a homomorphism of Lie groupoids over $M$.
We sometimes say for emphasis that $H$ is an {\em immersed}
subgroupoid of $G$. In case $\iota\!:H\ra G$ is an embedding,
we say that $H$ is an {\em embedded subgroupoid} of $G$.
A {\em closed subgroupoid} is a subgroupoid $H$ for which
$\iota\!:H\ra G$ is a closed embedding.
(In this paper, we shall only consider subgroupoids
over the same base.)

An immersed (source-connected) subgroupoid $H$ of $G$ is
completely determined by its Lie algebroid $\fh$, which is
a subalgebroid of $\fg$.
By Corollary \ref{cor:saif.8},
the corresponding foliation
$\cF(H)=\cF(\fh)$
of $G$ has trivial holonomy.
Conversely, we say that a Lie subalgebroid  $\fh$ of $\fg$
is integrable by a subgroupoid $H$ of $G$ if the associated
injective immersion $\iota\!:H\ra G$ integrates the inclusion
$\fh\ra\fg$. If such an integrating
subgroupoid exists, it is
unique.
We can therefore rephrase Corollary \ref{cor:saif.8}
as follows.

\begin{prop} \label{prop:sglg.3}
Let $G$ be a Lie groupoid with Lie algebroid $\fg$.
A subalgebroid $\fh\subset\fg$ can be integrated by a
subgroupoid of $G$ if and only if the foliation $\cF(\fh)$
has trivial holonomy.
\end{prop}

For an immersed subgroupoid $\iota\!:H\ra G$ and an arrow
$g\!:x\ra y$ in $G$, the (right) {\em coset}
$Hg$ is the immersed submanifold of $G$ given by
$$ H(y,\oo)\lra G\;,\;\;\;\;\;\;\;\;\;\;
   h\mapsto \iota(h)g\;.$$
These cosets are exactly the leaves of the associated
foliation $\cF(H)$ given by the Lie subalgebroid $\fh\subset\fg$
corresponding to $H$.
We will write $G/H=G/\cF(H)$
for the space of these right cosets,
with the quotient topology. Thus, $G/H$ is the quotient of $G$
obtained by identifying two arrows $g$ and $g'$ if and only if
$\src(g)=\src(g')$ and $g'g^{-1}$ belongs to (the image of) $H$.
The map $\src\!:G\ra M$ factors  as a map $G/H\ra M$,
and the right action of $G$
 on itself by multiplication induces a right
action on $G/H$ along this map $G/H\ra M$.

Recall that a foliation $\cF$ of a manifold $M$
is {\em simple} if it is given by the components of the
fibers of a submersion into a Hausdorff manifold.
It is called
{\em weakly simple} if there exists a smooth structure
of a possibly non-Hausdorff manifold on
$M/\cF$ such that the quotient map $M\ra M/\cF$ is
a submersion. A foliation $\cF$ of $M$
is {\em strictly simple}
if it is weakly simple and the space of leaves $M/\cF$
is Hausdorff.
In particular, any strictly simple foliation is simple,
and any simple foliation is weakly simple.
If the group of diffeomorphisms of $(M,\cF)$ acts
transitively on $M$, then
these three notions coincide
\cite[Theorem 4.3 (vi)]{MoerdijkMrcun2003}.

\begin{prop} \label{prop:sglg.5}
Let $H$ be a  subgroupoid
of a  Lie groupoid $G$,
and let $\cF(H)$ be the associated foliation of $G$.

(i) The subgroupoid $H$ is embedded in $G$
if and only if the foliation $\cF(H)$ is weakly simple.

(ii) The subgroupoid $H$ is closed in $G$
if and only if the foliation $\cF(H)$ is strictly simple.
\end{prop}

\begin{rem} \rm \label{rem:sglg.6}
In other words, the Lie subgroupoid
$H$ is embedded in $G$
if and only if there is a structure of a possibly
non-Hausdorff manifold on the
space of cosets $G/H$ such that the projection
$G\ra G/H$ is a submersion.
If this is the case, then
$H$ is closed if and only if $G/H$ is Hausdorff.
\end{rem}

\begin{proof}[Proof of Proposition \ref{prop:sglg.5}]
(i)
Suppose that $H$ is embedded, and consider the equivalence
relation $R$ on $G$ defining $G/H$. This relation is the
image of the map
$H\times_{M}G\ra G\times G$ sending a pair $(h,g)$ with
$\src(h)=\trg(g)$ to $(hg,g)$. Since $H$ is embedded
in $G$, $H\times G$ is also embedded in $G\times G$.
But $H\times_{M}G$ is a closed submanifold of $H\times G$,
therefore it follows that it is also embedded in $G\times G$.
Thus $R$ is an embedded submanifold of $G\times G$.
Both the projection $H\times_{M}G\ra G$ and the
composition $H\times_{M}G\ra G$ are submersions:
the first because it is a pull-back of the
submersion $\src\!:H\ra M$, the other because
it is isomorphic to the first.
By the Godement criterion \cite{Serre1992}
it follows that $G/H$ is a (possibly non-Hausdorff)
manifold such that $G\ra G/H$ is a submersion.

Conversely, suppose that the foliation $\cF(H)$ is
weakly simple, so $G/H$ has a structure of a possibly
non-Hausdorff manifold such that the quotient projection
$f\!:G\ra G/H$ is a submersion.
Since $\cF(H)$ refines the foliation of $G$ by the source-fibers,
the source map of $G$ factors through $f$ as a submersion
$\bar{\src}\!:G/H\ra M$.
The unit section $\uni\!:M\ra G$ induces a section
(hence an embedding) 
$f\com\uni\!:M\ra G/H$ of $\bar{\src}$.
Since the  groupoid $H$ fits into the pull-back square
$$
\xymatrix{
H \ar[d] \ar[r]^{\src} & M \ar[d]^{f\com\uni} \\
G \ar[r]^{f} & G/H
}
$$
it follows that $H$ is embedded in $G$.

(ii)
This follows directly from (i) and the fact that
$H$ is closed if and only if $G/H$ is Hausdorff.
\end{proof}

Recall that a
vector field $Y$ on a manifold $M$ is called {\em projectable}
with respect to a foliation $\cF$ of $M$
if its local flow preserves the
foliation, or equivalently, if the Lie derivative of
$Y$ along any vector field tangent to $\cF$ is
again tangent to $\cF$.
Following Molino \cite{Molino1977,Molino1988},
a foliation $\cF$ of $M$ is 
{\em transversely complete} if
any tangent vector on $M$
can be extended to a complete projectable vector field on $M$.
By Molino's structure theorem \cite{Molino1977},
the closures of the leaves of a transversely complete
foliation $\cF$ of $M$ are the fibers of a submersion $M\ra W$,
which is in fact a locally
trivial fiber bundle of (Lie) foliations.
We shall use (only) this property
of transversely complete foliations in this paper.
For our purpose 
it is also relevant to note that
any transversely complete foliation
has trivial holonomy.
Examples of transversely complete foliations
include foliations given by the fibers of locally trivial fiber bundles,
transversely parallelizable foliations on compact manifolds
\cite{Conlon1974,Molino1988}
and Lie foliations on compact manifolds
\cite{Fedida1971}.

An (immersed) subgroupoid $H$ of a Hausdorff
Lie groupoid $G$ is said
to be {\em transversely complete}
if the foliation $\cF(H)$ of $G$ by cosets of $H$
is transversely complete.
Recall that any right invariant foliation $\cF\subset\cF(\src)$
of $G$, which is transversely complete,
has trivial holonomy; hence Lemma \ref{lem:saif.3}
and Proposition \ref{prop:sglg.3} imply that any such
foliation $\cF$ is the foliation $\cF(H)$ associated
to a transversely complete subgroupoid $H$ of $G$.

\begin{exs} \rm \label{exs:sglg.9}
(1)
Suppose that $G$ is a connected
Lie group, thus a Lie groupoid
over a one-point space $M=\{\mathrm{pt}\}$.
Then any connected subgroup $H\subset G$ is
transversely complete.

(2)
Let $M$ be a manifold equipped with a foliation $\cF$.
Let $G$ be the pair groupoid $M\times M$. Then $\cF$ can be
viewed as a subalgebroid $\fh$ of the Lie algebroid
$\fg=T(M)$ of $M\times M$, and the corresponding
foliation $\cF(\fh)$ is
$\cF\times 0\subset T(M)\times T(M)=T(M\times M)$.
If $\cF$ has trivial holonomy,
then $\Hol(M,\cF)$ is an immersed
subgroupoid of $M\times M$, which
is transversely complete
whenever $\cF$ is. Thus, our terminology extends the usual
one for foliations.
\end{exs}

Recall that a Lie groupoid $G$ over $M$
is {\em transitive} 
if the map $(\trg,\src)\!:G\ra M\times M$
is a surjective submersion \cite{Mackenzie1987,MoerdijkMrcun2003}.
Any transitive Lie groupoid is automatically Hausdorff.

\begin{prop} \label{prop:sglg.13}
Any transitive subgroupoid of any Lie groupoid is
transversely complete.
\end{prop}

\begin{proof}
Let $H$ be a transitive subgroupoid of a Lie groupoid
$G$ over $M$. First we show that in this case $G$ is 
transitive as well. To see this, recall from
\cite[Proposition 5.14]{MoerdijkMrcun2003}
that a Lie groupoid over $M$
is transitive if and only if the restriction of
the target map to a source-fiber is a surjective submersion
onto $M$. 
Take any arrow $g\!:x\ra y$ in $G$.
Transitivity of $H$ implies that there exists an arrow
$h\!:x\ra y$ in $H$, and that the derivative
$(d\trg)_{h}\!:T_{h}(H(x,\oo))\ra M$ is surjective.
Since the right translation
$$ R_{h^{-1}g}\!:H(x,\oo)\lra G(x,\oo) $$
preserves the target and maps $h$ to $g$,
it follows that $(d\trg)_{g}\!:T_{h}(G(x,\oo))\ra M$
is surjective as well.
Thus we proved that $G$ is transitive, and in particular
Hausdorff.

We will now show that
the foliation $\cF(H)$ of $G$ associated to $H$
is transversely complete. 
Transitivity of $G$ implies that any
$x_{0}\in M$ has a neighbourhood $U$ and a section
$\sigma\!:U\ra G$ of the source map
such that $\trg(\sigma(x))=x_{0}$ for any $x\in U$.
The right translation
by this section provides a local
trivialization
$\src^{-1}(x_{0})\times U\ra \src^{-1}(U)$
of the source map, which also respects the foliation
$\cF(H)$ because $\cF(H)$ is right invariant.
From this it follows that it is enough to show
that the restriction of $\cF(H)$ to a source-fiber
is transversely complete.
Thus, for any arrow $g\!:x_{0}\ra y$
of $G$ and any $u\in T_{g}(G(x_{0},\oo))$
we need to find
a complete source-vertical
projectable vector field on $G(x_{0},\oo)$ with
value $u$ at $g$. 
Because the right translation
$R_{g}\!:\src^{-1}(y)\ra\src^{-1}(x_{0})$
preserves the foliation and sends $1_{y}$ into $g$,
we can assume without loss of generality that
$g=1_{x_{0}}$.

Since $H$ is transitive, the restriction of the target
map to $\src^{-1}(x_{0})$ is a surjective submersion onto $M$,
and hence we can write
$$ u=v+w $$
for a vector $v$ tangent to $H$ and a vector $w$
with $d\trg(w)=0$.
We can extend $v$ to a section of the Lie algebroid $\fh$ of $H$,
and we can assume that this section has compact support.
This section can be uniquely extended to a right invariant
source-vertical vector field $X$ on $G$, which is
tangent to $\cF(H)$ and has value $v$ at $1_{x_{0}}$.
The vector field $X$ is complete because its
anchor is \cite[p. 264]{KumperaSpencer1972}.
Next, we can uniquely extend $w$ to a left invariant
vector field $Y$ on $G(x_{0},\oo)$.
The flow of this vector field
is given by the family of right translations
$R_{\mathrm{exp}(tw)}$, $t\in\RR$, where
$\mathrm{exp}(tw)$ is the one-parameter subgroup
of the Lie group $G_{x_{0}}$ corresponding to $w$.
In particular, this flow is globally defined and
preserves the right
invariant foliation $\cF(H)$, so $Y$ is complete and
projectable. Furthermore, the vector field $Y$
commute with the right invariant vector field $X$.
This means that $X+Y$ is a complete projectable
vector field on $G(x_{0},\oo)$ and has value
$u$ at $1_{x_{0}}$.
\end{proof}

\begin{ex} \rm \label{ex:sglg.15}
Unlike the case where $G$ is a Lie group, 
an embedded subgroupoid of a Lie groupoid is not
necessarily closed, and its closure
may not be a subgroupoid. Indeed,
Proposition \ref{prop:sglg.5} implies that
the holonomy
groupoid of a simple, but not
strictly simple, foliation $\cF$ of a manifold $M$
is an embedded subgroupoid of the pair
groupoid $M\times M$ which is not closed.
If we take
$M=\RR^{2}\setminus (\{0\}\times [0,\infty))$
and if $\cF$ is the simple foliation of $M$
given by the second projection, the holonomy
groupoid of $(M,\cF)$ is embedded in $M\times M$,
but its closure is not a Lie subgroupoid.
\end{ex}

The following theorem shows
that transversely complete subgroupoids
behave very much like subgroups of
Lie groups.

\begin{theo} \label{theo:sglg.19}
Let $G$ be a  Hausdorff Lie groupoid,
and let $H$ be a transversely complete
subgroupoid of $G$.

(i)
The closure $\bar{H}$ of $H$ in $G$ is an embedded subgroupoid
of $G$.
    
(ii)
If $H$ itself is embedded in $G$, then it is closed, and
$G\ra G/H$ is a locally trivial fiber bundle.     
\end{theo}

\begin{proof}
(i)
By the Molino structure theorem referred to above,
the closures of the leaves of $\cF(H)$ are the fibers of
a locally trivial fiber bundle $\pi\!:G\ra W$.
Let $\bar{\src}\!:W\ra M$ be the submersion induced by
$\src$, with the section $\bar{\uni}=\pi\com\uni$.
Consider the following diagram:
$$
\xymatrix{
& P \ar[d]^{j} \ar[r] & M \ar[d]^{\bar{\uni}} \\
H \ar[r]^{\iota} \ar[ru] & G \ar[r]^{\pi} \ar[d]_{\src} & W \ar[dl]^{\bar{\src}} \\
& M &
}
$$
Here $P$ denotes the pull-back $\pi^{-1}(\bar{\uni}(M))$,
so $j\!:P\ra G$ is again a closed embedding.
Note that a point $g\!:x\ra y$ of $G$ belongs to $P$ precisely
when $g\in\bar{L}_{1_{x}}$, the closure of the leaf of $\cF(H)$
through the unit $1_{x}$. It follows that $P$ is a closed subgroupoid
of $G$. (If $g\in\bar{L}_{1_{x}}$ and $h\!:y\ra z$ belongs to 
$\bar{L}_{1_{y}}$, then also $hg\in\bar{L}_{g}$ by
right invariance of the foliation. Hence $hg\in\bar{L}_{1_{x}}$
because $\bar{L}_{g}=\bar{L}_{1_{x}}$.)

Also, the immersion $\iota\!:H\ra G$ evidently factors through $P$.
Finally, $H\ra P$ is dense, because the restriction of $H\ra P$ to
the source-fiber over a point $x\in M$ is the inclusion
$L_{1_{x}}\ra \bar{L}_{1_{x}}$, hence dense. In particular we have
$P=\bar{H}$.

(ii)
Note that the immersion of each leaf of $\cF(H)$
into the corresponding fiber of $\pi$ is
an embedding as well as dense. Thus the leaves of $\cF(H)$
are equal to the fibers of $\pi$,
and $H=P$ in the diagram above. It is clear that
$G\ra G/H$ is a locally trivial bundle in this case.
\end{proof}

\begin{ex} \rm \label{exs:sglg.39}
Let $G$ be a simply connected Lie group, and let
$\fg$ be the associated Lie algebra of right invariant
vector fields on $G$.
Suppose that $\cF$ is a Lie foliation on a manifold $M$,
given by the kernel of a Maurer-Cartan form
$$ \omega\in\Omega^{1}(M,\fg) $$
with values
in a $\fg$ \cite{Fedida1971,Molino1988}. 
This means that $\omega$ is non-singular and has
trivial formal curvature, $d\omega+\frac{1}{2}[\omega,\omega]=0$.
The foliation $\cF$ is
transversely parallelizable \cite{Conlon1974,Molino1988};
in fact, a vector field $Y$ on $M$ is projectable
with respect to $\cF$ if and only if the function $\omega(Y)$
is basic,  i.e.\ constant along the leaves of $\cF$.

The form $-\omega$ can be extended uniquely to
a flat connection form $\eta$ on the trivial principal
$G$-bundle $M\times G$
over $M$. (Here we take $-\omega$ instead of $\omega$ because
we use the Lie algebra of right invariant vector fields on $G$
instead the usual left invariant ones.) The kernel
of $\eta$ therefore defines a foliation $\cG$ of $M\times G$,
$$ (M\times G,\cG) \;,$$
which is $G$-invariant and has trivial holonomy.
Any leaf of $\cG$ is a covering space of $M$
with respect to the projection to $M$,
and the lift of $\cF$ to such a leaf is a strictly simple
foliation given by the fibers of the projection to $G$.
Any vector tangent to the foliation $\cG$
can be extended to a vector field tangent to $\cG$ of
the form $(X,\omega(X))$,
where $X$ is a (projectable) vector field on
$M$ for which $\omega(X)$ is a constant function
(and its value viewed as a right invariant vector field on $G$).

Consider the Lie groupoid
$$ M\times G\times M $$
over $M$,
i.e.\ the groupoid induced by $G$ along the
map $M\ra \{\mathrm{pt}\}$.
The arrows $x\ra y$ in $M\times G\times M$
are the triples $(y,g,x)$, where $g\in G$.
The 
foliation $\cG\times 0$ of $M\times G\times M$ is
right invariant
because $\cG$ is right $G$-invariant.
Thus $\cG\times 0$ is the foliation $\cF(\fh)$
of a subalgebroid $\fh$ of the Lie algebroid of
$M\times G\times M$,
$$ \cG\times 0=\cF(\fh)\;.$$
Furthermore, this foliation has trivial holonomy
because $\cG$ has trivial holonomy,
therefore it defines
a subgroupoid $H$ of the Lie groupoid
$M\times G\times M$ which integrates
Lie subalgebroid $\fh$.

Denote by $\Pi(M)$ the fundamental Lie groupoid of $M$,
and define a homomorphism of Lie groupoids
$$ \phi\!:\Pi(M)\lra M\times G\times M $$
over $M$ by path-lifting inside the leaves of $\cG$:
For any homotopy class $[\sigma]$ of a path $\sigma$
in $M$ from $x$ to $y$, let $\phi([\sigma])$ be the
unique lift of $\sigma$ which starts at
$(x,1,x)$ and lies inside a leaf of $\cG\times 0$.
The image of $\phi$ is precisely the groupoid $H$.
Using this path-lifting construction,
we see that an arrow $(y,g,x)$ of $M\times G\times M$
belongs to $H$ if and only if the points
$(y,g)$ and $(x,1)$ of $M\times G$ belong to the
same leaf of $\cG$.
Notice that $H$ is a transitive groupoid
and therefore transversely complete by
Proposition \ref{prop:sglg.13}.
Its isotropy group
$H_{x}$
at $x\in M$ is
the holonomy group of the Maurer-Cartan form
$\omega$.

By Theorem \ref{theo:sglg.19} (i) we know that
the closure
$\bar{H}$
of $H$ in $M\times G\times M$
is an embedded subgroupoid of $M\times G\times M$.
It is possible to give the following explicit
description of $\bar{H}$ in terms of the closures
$K_{x}$ of $H_{x}$ inside $G$ (such a description in fact
applies to any transitive subgroupoid):
An arrow $(y,g,x)\!:x\ra y$ in $M\times G\times M$
belongs to $\bar{H}$ if and only if it can be factored
as $g=hk$ for some $(y,h,x)\in H$ and $k\in K_{x}$.
In other words, the groupoid $H$ acts by conjugation
on the bundle $K$ of groups $K_{x}$ over $M$,
and $\bar{H}$ is constructed as the twisted
product of $H$ and $K$.
This example relates to the Cartier construction
of the groupoid closure
\cite{Cartier2004} in the context of Galois theory
for differential equations.
\end{ex}

\end{document}